\documentclass{elsart}
\usepackage{amsfonts,amsmath,amssymb}
\usepackage[francais,english]{babel}
\usepackage[dvips]{graphicx}
\usepackage[T1]{fontenc}
\begin{document}
\title{Some asymptotic results on density estimators by wavelet projections}
\author{Davit Varron, Université Catholique de Louvain}
\footnote{Institute of Statistics, voie du Roman Pays 20, 1348
Louvain la Neve, Belgium \textit{varron@stat.ucl.ac.be}}

\newcommand{\fn}{\hat{f}_n}
\newcommand{\und}{\vec{1}_d}
\newcommand{\suite}{_{n\geq 1}}
\newcommand{\ls}{\limsup}
\newcommand{\li}{\liminf}
\renewcommand{\e}{\epsilon}
\newcommand{\eni}{\po \e_{n,i}\pf_{n\geq1, i \le p_n}}
\newcommand{\tab}{_{n\geq1,\;i\le p_n}}
\newcommand{\lb}{\newline}
\newcommand{\dd}{\delta}
\newcommand{\ld}{\llcorner}
\newcommand{\lr}{\lrcorner}
\newcommand{\alp}{\alpha}
\newcommand{\lab}{\lambda}
\newcommand{\Lab}{\Lambda}
\newcommand{\gam}{\gamma}
\newcommand{\Gam}{\Gamma}
\newcommand{\sig}{\sigma}
\newcommand{\Sig}{\Sigma}
\newcommand{\mmi}{\mid\mid}
\newcommand{\mmmi}{\mid\mid\mid}
\newcommand{\mmii}{\mid\mid_\infty}
\newcommand{\gMid}{\Bigg |}
\newcommand{\Mid}{\Big |}
\newcommand{\Mmi}{\Big|\Big|}
\newcommand{\Mmii}{\Big|\Big|_\infty}
\newcommand{\T}{\Theta}
\newcommand{\tht}{\theta}
\newcommand{\sli}{\sum\limits}
\newcommand{\sliin}{\sum\limits_{i=1}^n}
\newcommand{\sliid}{\sum\limits_{i=1}^d}
\newcommand{\sliik}{\sum\limits_{i=1}^k}
\newcommand{\sliiN}{\sum\limits_{i=1}^N}
\newcommand{\slijk}{\sum\limits_{j=1}^k}
\newcommand{\proliin}{\prod\limits_{i=1}^n}
\newcommand{\proliid}{\prod\limits_{i=1}^d}
\newcommand{\proliik}{\prod\limits_{i=1}^k}
\newcommand{\proliiN}{\prod\limits_{i=1}^N}
\newcommand{\prolijk}{\prod\limits_{j=1}^k}
\newcommand{\ili}{\int\limits}
\newcommand{\proli}{\prod\limits}
\newcommand{\limn}{\lim_{n\rightarrow\infty}\;}
\newcommand{\limk}{\lim_{k\rightarrow\infty}\;}
\newcommand{\inv}{\frac{1}}
\newcommand{\lsn}{\limsup_{n\rightarrow\infty}}
\newcommand{\lin}{\liminf_{n\rightarrow\infty}}
\newcommand{\lsk}{\limsup_{k\rightarrow\infty}}
\newcommand{\lik}{\liminf_{k\rightarrow\infty}}
\newcommand{\bculi}{\bigcup\limits}
\newcommand{\bcali}{\bigcap\limits}
\newcommand{\wt}{\widetilde}
\newcommand{\indep}{{\bot}\kern-0.9em{\bot}}
\newcommand{\beq}{\begin{equation} }

\newcommand{\eeq}{\end{equation} }
\newcommand{\mcal}{\mathcal}
\newcommand{\sq}{\sqrt}
\newcommand{\rar}{\rightarrow}
\newcommand{\lar}{\leftarrow}
\newcommand{\cvps}{\rightarrow_{p.s.}\;}
\newcommand{\cvpr}{\rightarrow_{P}\;}
\newcommand{\cvloi}{\rightarrow_{\mcal{L}}\;}
\newcommand{\ii}{\vec{i}}
\newcommand{\jj}{\vec{j}}
\newcommand{\kk}{\vec{k}}
\newcommand{\zzi}{z_{\vec{i}}}
\newcommand{\zzij}{z_{\vec{i},\vec{j}}}
\newcommand{\zzI}{\vec{z}_I}
\newcommand{\zzIJ}{\vec{z}_{I,J}}
\newcommand{\zz}{z}
\newcommand{\nk}{n_k}
\newcommand{\nkm}{n_{k-1}}
\newcommand{\nkk}{n_{k+1}}
\newcommand{\nkd}{n_k}
\newcommand{\nkkd}{n_{k+1}}
\newcommand{\Ln}{L_{n,1}}
\newcommand{\Lnn}{L_{n,2}}
\newcommand{\Lnl}{L_{n,\ell}}
\newcommand{\Pn}{\Pi_{n,1}}
\newcommand{\Pnn}{\Pi_{n,2}}
\newcommand{\Pnl}{\Pi_{n,\ell}}
\newcommand{\mn}{\vec{m}_n}
\newcommand{\Mn}{\vec{M}_n}
\newcommand{\hhn}{h_{n,1}}
\newcommand{\hhnn}{h_{n,2}}
\newcommand{\hhk}{h_{n_k,1}}
\newcommand{\hhkk}{h_{n_k,2}}
\newcommand{\mk}{\vec{m}_k}
\newcommand{\Mk}{\vec{M}_k}
\newcommand{\hhnl}{h_{n,\ell}}
\newcommand{\hnk}{h_{n_k}}
\newcommand{\bn}{b_{n,1}}
\newcommand{\bnn}{b_{n,2}}
\newcommand{\bnl}{b_{n,\ell}}
\newcommand{\bnk}{b_{n_k}}
\newcommand{\bnkk}{b_{n_{k+1}}}
\newcommand{\hkl}{{h_{n_k,l}}}
\newcommand{\hnl}{{h_{n,l}}}
\newcommand{\Tn}{T_{n,1}}
\newcommand{\Tnn}{T_{n,2}}
\newcommand{\Tnl}{T_{n,\ll}}
\newcommand{\hamn}{\hat{m}_n}
\newcommand{\harn}{\hat{r}_n}
\newcommand{\hafn}{\hat{f}_n}
\newcommand{\srl}{\stackrel}
\newcommand{\wap}{(\Omega,\mathcal{A},\rm I\kern-2pt P)}
\newcommand{\aoo}{\Big\{}
\newcommand{\aff}{\Big\}}
\newcommand{\coo}{\Big [}
\newcommand{\cff}{\Big]}
\newcommand{\poo}{\Big (}
\newcommand{\pff}{\Big)}
\newcommand{\po}{\big (}
\renewcommand{\pf}{\big)}
\newcommand{\ao}{\big \{}
\newcommand{\af}{\big \}}
\newcommand{\co}{\big [}
\newcommand{\cf}{\big ]}
\newcommand{\pooo}{\bigg (}
\newcommand{\pfff}{\bigg)}
\newcommand{\aooo}{\bigg \{}
\newcommand{\afff}{\bigg \}}
\newcommand{\cooo}{\bigg [}
\newcommand{\cfff}{\bigg ]}
\newcommand{\poooo}{\Bigg (}
\newcommand{\pffff}{\Bigg)}
\newcommand{\aoooo}{\Bigg \{}
\newcommand{\affff}{\Bigg \}}
\newcommand{\coooo}{\Bigg [}
\newcommand{\cffff}{\Bigg ]}
\newcommand{\FF}{\mathcal{F}}
\newcommand{\TT}{\mathcal{T}}
\newcommand{\GG}{\mathcal{G}}
\newcommand{\BBGG}{\mathcal{B}(\mathcal{G})}
\newcommand{\CC}{\mathcal{C}}
\newcommand{\KK}{\mathcal{K}}
\newcommand{\SSS}{\mathcal{S}}
\newcommand{\BB}{\mathcal{B}}
\newcommand{\PP}{\mathcal{P}}
\newcommand{\HH}{\mathcal{H}}
\newcommand{\NN}{\mathcal{N}}
\newcommand{\MM}{\mathcal{M}}
\newcommand{\DD}{\mathcal{D}}
\newcommand{\cc}{\widetilde{c}}
\newcommand{\EEE}{\mathbb{E}}
\newcommand{\NNN}{\mathbb{N}}
\newcommand{\PPP}{ \mathbb{P}}
\newcommand{\CCC}{\mathbb{C}}
\newcommand{\KKK}{\mathbb{K}}
\newcommand{\RRR}{\mathbb{R}}
\newcommand{\ZZZ}{\mathbb{Z}}
\newcommand{\wtf}{\widetilde{f}}
\newcommand{\wtg}{\widetilde{g}}
\newcommand{\wth}{\mathfrak{h}}
\newcommand{\wtn}{\widetilde{n}}
\newcommand{\wtv}{\widetilde{v}}
\newcommand{\wtA}{\widetilde{A}}
\newcommand{\wtC}{\widetilde{C}}
\newcommand{\wtE}{\widetilde{E}}
\newcommand{\ovg}{\overline{g}}
\newcommand{\ovh}{\overline{h}}
\newcommand{\ovE}{\overline{E}}
\newcommand{\wtG}{\widetilde{G}}
\newcommand{\ovH}{\overline{H}}
\newcommand{\ovI}{\overline{I}}
\newcommand{\ovJ}{\overline{J}}
\newcommand{\ovK}{\overline{K}}
\newcommand{\wtN}{\widetilde{N}}
\newcommand{\wtP}{\widetilde{P}}
\newcommand{\ovR}{\overline{R}}
\newcommand{\wtPPP}{\widetilde{\mathbb{P}}}
\newcommand{\ovPPP}{\overline{\mathbb{P}}}
\newcommand{\wtF}{\widetilde{F}}
\newcommand{\wtI}{\widetilde{I}}
\newcommand{\wtK}{\widetilde{K}}
\newcommand{\wtFF}{\widetilde{\mathcal{F}}}
\newcommand{\ovFF}{\overline{\mathcal{F}}}
\newcommand{\cov}{\mathrm{Cov}}
\newcommand{\kif}{k\rightarrow\infty}
\newcommand{\nif}{n\rightarrow\infty}
\newcommand{\FFGG}{\FF\times\GG}
\newcommand{\vk}{\vskip10pt}
\newcommand{\Tproj}{\mathcal{T}_{0}}
\newcommand{\TTd}{\mathcal{T}_d}
\newcommand{\farc}{\frac}
\newcommand{\Nf}{\nabla_f}
\newcommand{\Nfn}{\nabla_f(\log_2n )}
\newcommand{\Cf}{\chi_f}
\newcommand{\nono}{\nonumber}
\newcommand{\Dn}{\Delta_n}
\newcommand{\DPn}{\Delta\Pi_n}
\newcommand{\Dalpn}{\Delta\alpha_n}
\newcommand{\DFn}{\Delta F_n}
\newcommand{\norm}{\mid\mid \cdot \mid\mid}
\newcommand{\Xt}{(X(t))_{t\in T}}
\newcommand{\zklj}{z_{k,l,j}}
\newcommand{\znlj}{z_{n,l,j}}
\newcommand{\Idd}{[0,1[^d}
\newcommand{\pipi}{p_{i,n,\ii_0}}
\newcommand{\ppd}{2^{-pd}}
\newcommand{\iideux}{\vec{i}\in {\{1,\ldots,2^p\}}^d}
\newcommand{\hk}{h_{n_k}}
\newcommand{\zik}{z_{i,n_{k}}}
\newtheorem{theo}{Theorem}
\newtheorem{ptheo}{Proof of theorem}
\newtheorem{plem}{Proof of lemma}[section]
\newtheorem{preuveprop}{Proof of property}[section]
\newtheorem{defi}{Definition}[section]
\newtheorem{propo}{Proposition}[section]
\newtheorem{popo}{Proof of proposition}[section]
\newtheorem{coro}{Corollary}[theo]
\newtheorem{pcoro}{Proof of corollary}[theo]
\newtheorem{Fact}{Fact}[section]
\numberwithin{equation}{section}

 \maketitle
\begin{center}\textbf{Abstract}\end{center}
\begin{small}Let $(X_i)_{i\geq 1}$ be an i.i.d. sample on $\RRR^d$ having
density $f$. Given a real function $\phi$ on $\RRR^d$ with finite
variation and given an integer valued sequence $(j_n)$, let $\fn$
denote the estimator of $f$ by wavelet projection based on $\phi$
and with multiresolution level equal to $j_n$. We provide exact
rates of almost sure convergence to 0 of the quantity $\sup_{x\in
H}\mid \fn(x)-\EEE(\fn)(x)\mid$, when $n2^{-dj_n}/\log n \rar
\infty$ and $H$ is a given hypercube of $\RRR^d$. We then show
that, if $n2^{-dj_n}/\log n \rar c$ for a constant $c>0$, then the
quantity $\sup_{x\in H}\mid \fn(x)-f\mid$ almost surely fails to
converge to 0. \lb
\end{small}
\section{Introduction and statement of the results}
The well known wavelet theory (see, e.g., Mallat (1989)) has
proven useful in may branches of applied mathematics and
functional estimation in the field of statitics. In this paper, we
are interested in estimating the Lebesgue density $f$ of an
independent, identically distributed (i.i.d.) sample $(X_i)_{i\geq
1}$ on $\RRR^d$. Let $\phi$ be a (mother wavelet) real function on
$\RRR^d$. There exists $N=2^d-1$ associated (father wavelet
functions) $\Psi_{j},\;1\le j\le N$, so as any function $f\in
L^2(\RRR^d,\lab)$ has the following orthogonal representation (for
more details, see, e.g., Masry (1997)). \begin{align}
f:=&\sli_{k\in \ZZZ^d} \alp_{j_0,k}\phi_{j_0,k}+\sli_{j\geq
j_0}\sli_{i=1}^N\sli_{k\in
\ZZZ^d} \beta_{i,j,k} \Psi_{i,j,k},\label{expan}\text{ where} \\
\phi_{j,k}:=&2^{dj/2}\phi(2^j\cdot-k),\;\Psi_{i,j,k}:=2^{dj/2}\Psi_i(2^j\cdot-k),\;
i\le N,\;j\geq 1,\;k\in \ZZZ^d.\end{align} \subsection{The linear
wavelet projection estimator} The linear wavelet projection
estimator (see, e.g., Masry (1997)) of $f$ is constructed by
estimating the coefficients $\alp_{j_0,k},\;\beta_{i,j,k}$ by
their empirical analogues: \beq
\hat{\alp}_{j_0,k}:=\frac{1}{n}\sliin
\phi_{j_0,k}(X_i),\;\hat{\beta}_{i,j,k}:=\frac{1}{n}\sliin
\Psi_{i,j,k}(X_i),\eeq and stopping the expansion (\ref{expan}) at
a deterministic (multiresolution) level $j_n$, which will be
assumed to grow with the sample size $n$. \beq \fn(x):=\sli_{k\in
\ZZZ^d}\hat{\alp}_{j_0}\phi_{j_0,k}+\sli_{j=j_0}^{j_n-1}\sli_{i=1}^N\sli_{k\in\ZZZ^d}\hat{\beta}_{i,j,k}\Psi_{i,j,k}.\eeq
The aim of this paper is to describe the almost sure asymptotic
behaviour of the quantity $\mid \fn(x)-f(x)\mid$, uniformly in
$x\in H$, where $H$ is a given an hypercube of $\RRR^d$.
Obviously, the asymptotic behaviour of $(j_n)\suite$ plays a
crucial role, and $2^{-d j_n}$ can be intuitively compared to the
bandwidth when estimating $f$ by usual kernel methods. Massiani
(2003) has given an asymptotic result of $\fn$ when the sample
$(X_i)_{i\geq 1}$ takes values in $\RRR$ and under the following
conditions, with $h_n:=2^{-d j_n}$ :
\begin{align}
& h_n\downarrow
0,\;\;\;nh_n\uparrow\infty,\;\;n h_n/\log n\rar
\infty,\;\log(1/h_n)/\log\log n\rar \infty,\label{CRS}\\
&f \;\mathrm{is\;continuous\; and\; strictly \; positive\;
on\;an\;open\;subset\;}O\;\mathrm{and}\;H\subset O,\label{regu}\\
&\phi\;\mathrm{has\;finite\;variation\;on}\;\RRR^d\mathrm{\;and\;has\;a\;compact\;support}.\label{phi1}
\end{align}
Conditions (\ref{CRS}) are called the Csörg\H{o}-Révész-Stute
conditions. Massiani proved that, under (\ref{CRS}), (\ref{regu})
and (\ref{phi1}) we have, almost surely, \beq \limn \sup_{x\in H}
\pm {\poo\frac{n 2^{d j_n}}{2 f(x)\log(2^{dj_n})}\pff}^{1/2}\poo
\fn(x)-\EEE\po \fn(x)\pf\pff=1.\label{resultat}\eeq We also refer
to Masry(1997) for related results when $(X_i)_{i\geq 1}$ is a
stationary strongly mixing sequence. To prove (\ref{resultat}),
the author made use of the following expression of $\fn$ (see,
e.g, Masry (1997)) \begin{align}\fn(x):= &\frac{1}{n}\sliin
K_{j_n}(x,X_i),\label{exprfn}\text{ where}\\
K(x,y):=&\sli_{k\in\ZZZ^d}\phi(x-k)\phi(y-k),\;x,y\in \RRR^d \label{K},\\
K_{j_n}(x,y):=&2^{dj_n}K\po 2^{j_n}x,2^{j_n}y\pf\label{Kjn}.
\end{align}
Then, the author showed that (\ref{exprfn}) can be expressed quite
simply  with the functional increments of the empirical
distribution function, and made extensively use of related results
established by Deheuvels and Mason (1992). We point out the fact
that the just mentioned pioneering results do not cover the case
where the sample is multivariate ($d>1$), as this result relies on
the strong approximation theorem of Koml\'os \textit{et al.}
(1977). As a consequence, Massiani could only prove
(\ref{resultat}) when $d=1$. However, Mason (2004) recently made a
skillful use of some recent tools in empirical processes theory to
extend the results of Deheuvels and Mason (1992) to a more general
framework, which covers the case where $d>1$. As a consequence, we
are now able to prove the following result.
\begin{theo}\label{TDA1}
Under assertions (\ref{CRS}), (\ref{regu}), (\ref{phi1}) we have
almost surely:
\begin{align}
\nono
(i)&\;\mathrm{For\;each\;}\e>0,\mathrm{\;there\;exists\;}n(\e)\mathrm{\;such\;that,\;for
\;each}\;n\geq n(\e)\;\mathrm{and\;}\;x\in H,\;\\
\nono &{\poo\frac{n 2^{d j_n}}{2 f(x)\log(2^{dj_n})}\pff}^{1/2}\poo
\fn(x)-\EEE\po \fn(x)\pf\pff \in [-1-\e,1+\e],\\
 \nono (ii)&\;
\mathrm{For\;each\;}v\in [-1,1]\mathrm{\; and\;}
\e>0,\mathrm{\;there\;exists\;}n(\e,v)\mathrm{\;such\;that,\;for
\;each\;}n\geq n(\e,v),\;\\
\nono &\inf_{x\in H}\Mid{\poo\frac{n 2^{d j_n}}{2
f(x)\log(2^{dj_n})}\pff}^{1/2}\poo \fn(x)-\EEE\po
\fn(x)\pf\pff-v\Mid<\e.
\end{align}
\end{theo}
As mentioned above, the uniform behaviour of the increments of the
empirical process shows up to rule that of $\fn(x)$. Moreover, it
is well known (see Deheuvels and Mason (1992)) that this behaviour
changes abruptly when conditions (\ref{CRS}) are replaced by the
following Erdös-Rényi conditions: \beq h_n\downarrow 0,\;nh_n
\uparrow \infty,\;nh_n/\log n\rar c.\label{ER}\eeq Here, $c>0$ is
a finite constant. Since the pioneering result of Deheuvels and
Mason (1992), several extensions have been made. In Varron (2007),
Varron recently showed that this nonstandard UFLL still holds when
$d>1$. Our next result shows that, under (\ref{ER}), the
nonstandard behaviour of the empirical increments implies that the
uniform strong consistency of $\fn$ on a hypercube $H$
\textit{fails to hold}.
\begin{theo}\label{TDA2}
Under (\ref{regu}), (\ref{phi1}), (\ref{ER}), the following event
has probability 1:
\begin{align}\exists\;\e>0,\;\forall n_0,\;\exists n\geq n_0,\;\exists\;x_n\in H\;\mathrm{fulfilling }\;\Mid \frac{\fn(x_n)}{f(x_n)}-1\Mid>\e\label{nonconsistance}.
\end{align}
\end{theo}
\section{Proofs of Theorem \ref{TDA1}}\label{deux}
Recall that $h_n:=2^{-dj_n},\;n\geq 1$. To prove Theorem
\ref{TDA1}, we shall require some more notations. Given
$s:=(s_1,\ldots,s_d)$ and $v:=(v_1,\ldots,v_d)$, we shall write
$s\prec v$ whenever $s_i\le v_i$ for each $1\le i\le d$ and we
shall write $[s,v]$ for the set $\{u\in \RRR^d,\;s\prec u\prec
v\}$. The increments of the empirical process based on
$(X_i)_{i\geq 1}$ are defined as follows ($C$ denoting a Borel
subset of $\RRR^d$)\beq \Dalpn(x,h_n,C):= n^{1/2}\pooo
\frac{1}{n}\sliin 1_{C}\poo
\frac{X_i-z}{h_n^{1/d}}\pff-\EEE\poo1_{C}\poo
\frac{X_i-z}{h_n^{1/d}}\pff\pff\pfff.\label{Dalpn}\eeq A standard
argument of  homothety shows that we can make the following
assumption with no loss of generality: \beq
\phi\;\mathrm{has\;finite\;variation\;on}\;\RRR^d\mathrm{\;and\;has\;a\;support\;included\;in
}\;[-1/4,1/4]^d.\label{phi2}\eeq Set $I^d:=[-1/2,1/2]^d$, and
consider the space $B(I^d)$ of real, bounded, Borel functions on
$I^d$. We endow $B(I^d)$ with the usual supremum norm, namely
$\mmi g\mmi:=\sup\{\mid g(s)\mid,\;s\in I^d\}$. The proof of
Theorem \ref{TDA1} strongly relies on the following fact, which is
due to Mason (2004). Call $H_n$ the set of points $x\in H$ such
that $2^{j_n}x\in \ZZZ^d$, and define the following Strassen-type
set: \begin{align} \nono\SSS_{I^d}:=&\aoo g\in
B(I^d),\;\exists\dot g\;\mathrm{Borel}, \; \ili_{I^d} {\dot g}^2
d\lab\le 1,\;\forall s \in I^d,\; g(s)=\ili_{[s,1/2]} \dot g
d\lab\aff.\label{strassen}\end{align}  \textbf{Fact 1 (Mason,
2004)}: Set \beq g_{n,x}(s):=\frac{\Dalpn(x,h_n,[s,1/2])}{\sqrt{2
f(x)h_n \log(1/h_n)}},\;x\in H,\;n\geq 1,\;s\in
I^d.\label{gnx}\eeq Under assumption (\ref{CRS}) and (\ref{regu}),
we have almost surely:
\begin{align}
\nono (a)&\;\forall\e>0,\;\exists n(\e),\;\forall n\geq n(\e)\;\mathrm{and\;}\;x\in H,\;\inf_{g\in \SSS_{I^d}}\mmi g_{n,x}-g\mmi<\e\\
 \nono (b)&\;
\forall g\in \SSS_{I^d}\mathrm{\; and\;} \e>0,\;\exists
n(\e,g),\;\forall n\geq n(\e,g),\;\inf_{x\in H_n} \mmi
g_{n,x}-g\mmi<\e.
\end{align}
This fact is a nearly direct consequence of Theorem 1 of Mason
(2004), considering the class $\FF:=\{1_{[s,1/2]},\;s\in I^d\}$,
by Remark $F.2$ in Mason (2004).\lb
\begin{rem}
We point out the fact that Theorem 1 in Mason (2004) cannot lead
to Fact 1 directly, because $(b)$ involves the quantity
$\inf\{\mmi g_{n,x}-g\mmi,\;x\in H_n\}$ instead of $\inf\{\mmi
g_{n,x}-g\mmi,\;x\in H\}$. However, looking closely at the proof
of point $(b)$ of Theorem 1 in Mason (2004), we can see that $H$
can be replaced by $H_n$, as we can choose
$\{z_{1,n},\ldots,z_{m_n,n}\}:=H_n$ in his proof of Lemma 2.
\end{rem}
Set, for fixed $x\in \RRR^d$ and $n\geq 1$ (recall (\ref{K})),
\begin{align}
\wtK_{n,x}(s):=&K(2^{j_n}x,2^{j_n}x+s),\; n\in \RRR^d,\;s\in \RRR^d\label{wtKnx}\\
\sig^2_{n,x}:=&\ili_{\RRR^d} {\wtK}^2_{n,x}(s) ds\label{signx}
\end{align}
By assumption (\ref{phi2}), each $\wtK_{n,x}$ has support included
in $I^d$. Now, we consider the following continuous linear
applications, from $\po B(I^d),\norm\pf$ to $\RRR$. For fixed
$x\in \RRR^d$ and $n\geq 1$, set \beq \T_{n,x}\po
g\pf:={\sig_{n,x}}^{-1}\ili_{I^d}g(s) d\wtK_{n,x}(s),\;g\in
B(I^d).\label{Tnx}\eeq With these notations, we obviously have for
each $n\geq 1$ and $x\in\RRR^d$, almost surely, \beq
{\poo\frac{nh_n}{2f(x)\log(1/h_n)}\pff}^{1/2} \poo
\fn(x)-\EEE(\fn(x))\pff= \T_{n,x}(g_{n,x})\label{relation},\eeq so
as the random objects involved in Theorem \ref{TDA1} show up to be
correctly chosen functions of the increments of the empirical
process.\lb We first focus on proving point $(i)$ of Theorem
\ref{TDA1}. Standard analysis shows that \beq
\T_{n,x}\po\SSS_{I^d}\pf= [-1,1]\text{  for each $x\in \RRR^d$ and
$n\geq 1$}.\label{inclus}\eeq Moreover, by definition of the
$\wtK_{n,x}$ and by (\ref{phi2}) we have, $\Upsilon$ denoting the
total variation of a function,\begin{align} &\sup_{n\geq 1,\;x\in
H}
\Upsilon\po\wtK_{n,x}\pf<\infty\label{Var},&\\
&\inf_{n\geq 1,\;x\in H}
\sig_{n,x}>0,\;&\mathrm{whence}\label{integ}\\
&\sup_{n\geq 1,\;x\in H}\;\sup_{g\in B(I^d),\;\mmi g\mmi =1} \mid
\T_{n,x}(g)\mid<\infty\label{lip}.&
\end{align}
Note that (\ref{integ}) is a consequence of the Cauchy-Schwartz
inequality, as $\ili\wtK_{n,x}(s)ds=1$ (see, e.g., Meyer (1990),
p. 33). Now, combining (\ref{lip}) and point $(a)$ of Fact 1, we
conclude that point $(a)$ of Theorem \ref{TDA1} is true, by
routine topology.\lb We shall now prove point $(b)$ of Theorem
\ref{TDA1}. Recall that $x\in H_n$ if an only if $2^{j_n}x\in
\ZZZ^d$. Hence, by definitions (\ref{K}) and (\ref{wtKnx}) we have
\beq\wtK_{n,x}=\wtK_{0,0},\;\;\T_{n,x}=\T_{0,0},\;n\geq 1,\;x\in
H_n\label{egalite}.\eeq Now fix $\e>0$ and $v\in [-1,1]$.
Recalling (\ref{inclus}) we choose $g\in \SSS_{I^d}$ fulfilling
$\T_{0,0}(g)=v$. Now, as $\T_{0,0}$ is Lipschitz, and by point
$(b)$ of Fact 1, we conclude that, almost surely, there exists
$n(\e,g)$ such that, for each $n\geq n(\e,g)$, there exists
$x_n\in H_n$ fulfilling
$\mid\T_{n,x}(g_{n,x})-v\mid=\mid\T_{0,0}(g_{n,x})-\T_{0,0}(g)\mid<\e.$
The end of the proof follows readily, as $[-1,1]$ is compact.
$\Box$
\section{Proof of
Theorem \ref{TDA2}}\label{trois} In this section, condition
(\ref{CRS}) is replaced by condition (\ref{ER}). We first define
\begin{align} \DFn(x,h_n,C):=&\frac{1}{c
f(x)n h_n}\sliin
1_{C}\poo\frac{X_i-x}{h_n^{1/d}}\pff,\;C\;\mathrm{Borel},\;x\in
H,\;n\geq 1,\label{DFn}\\
\wtg_{n,x}(s):=&\DFn(x,h_n,[s,1/2]),\;s\in I^d,\;x\in H,\;n\geq
1.\label{wtgnx}
\end{align}
Set $ h(x)=
    (x\log x-x+1)1_{(0,\infty)}(x)+1_{\{0\}}(x)$
 for $x\geq 0$ and $h(x)=\infty$ otherwise. Now consider the following limit sets depending on a real
parameter $v>0$ : \begin{align} \Gam_{v,I^d}:=&\aoo
g,\;\exists\dot g\;\mathrm{Borel},\;\ili_{I^d} h\po{\dot g}\pf
d\lab\le 1/v,\;\forall s \in I^d,\; g(s)=\ili_{[s,1/2]} \dot g
d\lab\aff.\label{Gamv}\end{align} We shall make use of the
following result, which is a consequence of Theorem 1 of Varron
(2007). Recall that $x\in H_n$ if and only if $2^{j_n}x\in
\ZZZ^d$.\lb \textbf{Fact 2 (Varron)} Under assumptions
(\ref{regu}) and (\ref{ER}), the following assertions hold with
probability one.
\begin{align}
\nono(a)&\;\forall \e>0,\;\exists n(\e),\;
\forall n\geq n(\e)\;\mathrm{and \;}x\in H,\;\inf_{g\in \Gam_{cf(x),I^d}}\mmi \wtg_{n,x}-g\mmi<\e;\\
\nono (b)&\; \forall\;x\in H\;\mathrm{,}\;g\in \Gam_{cf(x)},\;
\e>0,\;\exists n(\e,g),\;\forall n\geq n(\e,g),\; \inf_{x\in H_n}
\mmi \wtg_{n,x}-g\mmi<\e.
\end{align}
\begin{rem}
Note that Fact 2 differs from Theorem 1 in Varron (2007) by two
aspects. First, the involved class of set is
$\FF_1:=\{1_{[s,1/2]},\;s\in I^d\}$ instead of
$\FF_2:=\{1_{[0,s]},\;s\in [0,1]^d\}$. However, by a standard
translation argument, one can trivially transpose Theorem 1 in
Varron (2007) from $\FF_2$ to $\FF_1$. Second, the cube $H$ is
replaced by $H_n$ in point $(b)$. As in Remark 1, we underline
that this replacement can be made by a close look at the arguments
of Varron (2007)
\end{rem}
To prove Theorem \ref{TDA2}, we shall make use of point $(b)$ of
Fact 2. Similarly to what was done in \S \ref{deux}, we introduce
the following linear applications \beq
\T'_{n,x}(g):=\ili_{I^d}g(s) d\wtK_{n,x}(s),\;g\in
B(I^d)\label{Tpnx}.\eeq Notice that, for any $n\geq 1$ and $x\in
H_n$ we have $\T'_{n,x}=\T'_{0,0}$ since $\wtK_{n,x}=\wtK_{0,0}$.
Now, as \beq
\frac{\fn(x)}{f(x)}:=\T'_{n,x}(\wtg_{n,x})=\T'_{0,0}(\wtg_{n,x}),\;n\geq
1,\;x\in H_n,\label{ex2}\eeq the proof of point $(b)$ of Theorem
\ref{TDA2} would be a direct consequence of point $(b)$ of Fact 2,
provided that the following statement is true for some $\dd>0$ :
\beq \bcali_{x\in
H}\T'_{0,0}\po\Gam_{cf(x),I^d}\pf\supset[1-\dd,1+\dd].\eeq Now, by
definition of $\Gam_{v,I^d},\;v>0$ we obviously have
$$\bcali_{x\in H}\T'_{0,0}\po\Gam_{cf(x),I^d}\pf =
\T'_{0,0}(\Gam_{c f(x_0),I^d})=:J,$$ where
$f(x_0)=\sup\{f(x),\;x\in H\}$. Note that, when $d=1$, the set $J$
can be described by making use of the optimisation techniques of
Deheuvels and Mason (see Deheuvels and Mason (1991), Theorem 3 and
4 and Deheuvels and Mason (1992), Theorem 4.2). To conclude the
proof of Theorem \ref{TDA2}, we shall now show that $J$ has a
nonempty interior. Define the following function for $I^d$ to
$\RRR$: \beq
  g_0:\;(s_1,\ldots,s_d)  \rar  \proliid(\frac{1}{2}-s_i). \eeq
Obviously, $g_0$ belongs to $\Gam_{c f(x_0),I^d}$, as $\dot
g_0\equiv 1$ fulfills the requirements stated in (\ref{Gamv}).
Moreover, an integration by parts leads to the conclusion that
\beq T'_{0,0}(g_0)=\ili_{\RRR^d}
\wt{K}_{n,x}(s)ds=\ili_{\RRR^d}K(s)ds=1.\label{rsl}\eeq As
$\Gam_{cf(x_0),I^d}$ is convex and $T'_{0,0}$ is linear, the set
$J$ is an interval that contains $T'_{0,0}(g_0)=1.$ Moreover, as
$h$ is continuous at $x=1$, we have $\rho g_0 \in \Gam_{c
f(x_0),I^d}$ and $\rho^{-1} g_0 \in \Gam_{c f(x_0),I^d}$ for
$\rho>1$ small enough, which entails, by linearity of $\T'_{0,0}$,
\beq \inf_{g\in\Gam_{c f(x_0),I^d}} \T_{0,0}(g)\le
\rho^{-1}<1<\rho\le\sup_{g\in\Gam_{c f(x_0),I^d}}
\T_{0,0}(g)\Box\eeq

\end{document}